\magnification=1200
\null\bigskip\bigskip\bigskip
\centerline{THE NATURE OF PROPERLY HUMAN MATHEMATICS.}
\bigskip
\centerline{David Ruelle.\footnote{$\dagger$}{IHES, 91440 Bures sur Yvette, France. email: ruelle@ihes.fr}.}
\bigskip\bigskip
\centerline{This note is dedicated to the memory of Jean Ginibre who}

\centerline{expressed his interest for structures in his publications [8].}
\bigskip\bigskip
	{\it Abstract.}  We claim that human mathematics is only a limited part of the consequences of the chosen basic axioms.  Properly human mathematics varies with time but appears to have universal features which we try to analyse.  In particular the functioning of the human brain privileges concept naming and short formulations.  This leads to organising mathematical knowledge structurally.  We consider briefly the problem of non-mathematical sciences.
\vfill\eject
\bigskip\bigskip
	The functioning of the human brain is very different from that of electronic computers.  Various mathematicians have noted that in doing mathematics the brain uses natural languages, functions combinatorially, is slow, prone to error, has very limited working memory, is very highly parallel, and uses unconscious thinking.  Inherent to mathematics is the fact that there are some very long proofs of theorems with short formulation.  We argue that the way the human brain functions explains some features of human mathematics, as opposed to listing general consequences of the axioms.  In particular, the human brain privileges concept naming and short formulations.  This leads to a structural organisation of mathematical knowledge.
\medskip\noindent
{\bf 1. A non-philosophical approach.}
\medskip
	It is natural for scientists to reflect on their handling of scientific problems.  An example of this is the article [6] by Jaffe and Quinn.  The following note is a modest attempt in this direction (after earlier small contributions, for instance [9]).  We are concerned here with finding natural limitations to human scientific progress and specifically in human mathematics.  Studying such limitations involves estimating both human intellectual capabilities and the complexity of things investigated.
\medskip
	We are addressing here a philosophical problem, and we want to make clear why we want to stay away from the traditions of philosophical literature.  A couple of citations by scientists may illustrate the problem.  There is the popular approximate quote attributed to C.F. Gauss:

\qquad When a philosopher says something that is true then it is trivial. When he says

\qquad something that is not trivial then it is false.

\noindent
A.N. Whitehead also wrote [14]:

\qquad The safest general characterisation of the European philosophical tradition is

\qquad that it consists of a series of footnotes to Plato.

\noindent
My view is that philosophical literature uses concepts like `understanding', `thinking', `meaning', `existence', `consciousness', `God', which refer to an introspective assessment of the functioning of the human brain.  I prefer to make my own assessment, staying close to facts and away from philosophical tradition.
\medskip
	Note that `there exists' is a well-defined mathematical concept to be used in a proper setting, we shall also try to make sense in Section 4 of `understanding a mathematical proof'\footnote{*}{This is not the place to discuss delicate questions, like G\"odel's ``proof of the existence of God''.  See Wikipedia for reference and discussion; this ``proof'' was pointed out to me by C. Liverani.}.  In practice many references to `understanding', `meaning', etc., are without ill consequences (and hard to avoid).  I think however that it is important to notice that some arguments like {\it cogito ergo sum} are based on introspection and to be treated as psychological evidence.  Such arguments are not to be used as logical proof.
\medskip\noindent
{\bf 2. What is human mathematics?}
\medskip
	Mathematics consists in deriving consequences (theorems) from a set of assumptions (axioms) by application of given logical rules.  The set of axioms mostly used currently is ZFC (Zermelo-Fraenkel-Choice set-theoretical axioms).  Axioms and theorems can be formulated in a formal language.  ZFC is fairly believable by mathematicians (a typical axiom is `there exists an infinite set').  We remind the reader that the consistency of ZFC cannot be proved (this follows from G\"odel's incompleteness theorems).
\medskip
	Human mathematics is based on natural languages (ancient Greek, English, etc.) which can in principle be translated into formal language (but is hardly understandable after translation).  There also exist {\it formal proofs} [3],[4] using a computer to verify the correctness of deductive steps created by human mathematicians (these mathematicians use a {\it proof assistant} like {\it Coq} or {\it HOL Light}).  Formal proofs provide logical certainty if the basic axioms are consistent and the computer functioning is reliable.  This is an advantage over the usual human proofs using a natural language since long proofs often contain logical errors or confusing statements due to the informality of natural languages.
\medskip
	The advantages of formal proofs, or other proofs using computers in an essential way, come at a the cost of {\it human understandability}.  Now understanding\footnote{*}{M.A.F. Sanju\'an [10] among others has discussed the status of `understandability'.} is one of the `philosophical' concepts which involves introspective reference to the functioning of the human brain.  We shall not try to give a `scientific' definition of understanding based on the functioning of the human brain.  We shall instead use generally accepted or acceptable features of mathematical understanding based on verifiable facts, see the example of `understanding a mathematical proof' in Section 4 below.  Note that a long computer proof may be verifiable step by step by a human, but this does not qualify as global human understanding.
\medskip
	Properly human proofs are thus using a natural human language, mathematical jargon, and formulae, not translated into formal language, and without essential use of a computer.  I am reminded here of a statement by P. Deligne that what he was interested in were mathematical proofs that he could understand, excluding computer proofs or proofs so long that he could not understand them globally.  [Quotation autorised by Deligne who adds the following:  ``My ideal is to be able to prove (to myself) everything I state, and at least I keep in mind the exceptions I made.  More important perhaps, the utility of a proof is double: ensure that a statement is true, and understand it (in particular by relating it to other statements).  Hence the adage that I learned from N. Katz: what is worth proving is worth proving again.'']
\medskip\noindent
{\bf 3. Features of human mathematical thinking.}
\medskip	
	Properly human mathematics {\it uses natural languages} (plus jargon and formulae) perhaps with some noncreative use of computers.  It is usually in written form, with the legacy of earlier times.
\medskip
	One general feature of mathematics, which applies in particular to properly human mathematics, is that {\it some theorems have very long proofs} with respect to the length of their statement.  This is related to G\"odel's incompleteness theorems; for instance one may invoke the fact that the halting problem for a Turing machine is undecidable [11].  A consequence of this is that the rate of progress per mathematician (using properly human mathematics) should be decreasing with time.  The rate of progress in other sciences is influenced by technological and political factors.
\medskip
	Comparison between human and computer intelligence have been made by Turing [12] and von Neumann [13].  Compared to a computer, {\it the brain is slow, prone to error, has limited memory, and is very highly parallel}.
\medskip
	[The reaction time of a neurone is between 10$^{-4}$ and 10$^{-2}$ seconds (von Neumann estimate) while current computer characteristic times are 10$^{-6}$ to 10$^{-9}$ seconds and faster.  There are different kinds of human memory, in part huge, but the short-term memory (working memory) involved in `thinking' is quite limited (to about 7 objets) compared with the huge short-term memory of computers.]
\medskip
	We have reflections of various mathematicians on mathematical creativity (Gauss, Poincar\'e [7], Hadamard [2], Hardy [5], etc.).  They stress the role of {\it unconscious thinking}, {\it limited working memory}, and the {\it combinatorial nature of mathematical thinking}.  These reflections do not amount to a theory of mathematical creativity (otherwise we could program a creative mathematical computer).
\medskip
	To summarize, we have noted that in properly human mathematics {\it the brain uses natural languages, functions combinatorially, is slow, prone to error, has very limited working memory, is very highly parallel, and uses unconscious thinking.  Furthermore there are some very long proofs}.
\medskip\noindent
{\bf 4. On the nature of human mathematical creativity: structures.}
\medskip
	We remark that human mathematics does not consist of a list of all logical consequences of the axioms, but is based on a structural understanding of some of these consequences.  The structures we refer to are things like geometric objects, groups, compact spaces, etc., that are generally not visible in the axioms.  A structural presentation of mathematics has for instance been advocated by Bourbaki (but this did not include categories and functors), Grothendieck, and Voevodsky.  The structural understanding of mathematics is a human creation which changes with time and depends in some way on the features of the human mind.
\medskip
	How do we account for the use of structures in human mathematics?  The use by the brain of a natural language, slowness, and limited working memory makes natural the use of short formulations (words, short sequences of symbols).  This, together with logical necessity, explains at least in part the use of concepts and structures in mathematical thinking.  
\medskip
	We see mathematical thinking as creating or checking a text in human mathematical language spoken or written in the natural direction.  Specifically we mean that certain sets of small size pieces extracted (according to suitable rules) from earlier parts of the text should satisfy suitable rules (including number, length and nature of the pieces).  This involves short-term and longer-term memories and makes sense of `understanding a mathematical proof'.  Except for dealing with the idiosyncrasies of human language, checking a human mathematical text could easily be done by computer (and something like this is essentially done in formal proofs).  Creating an interesting mathematical theorem is however at this time basically beyond the reach of computers.  In particular the role of high brain parallelism and the use of unconscious thinking remain unclear.  This is where the human versus computer problem [12] stands at this moment for mathematics.
\medskip
	Apart from the use of natural languages, short names for structures, and combinatorial functioning, human mathematics has to deal with very long proofs.  The proof [1] of the odd order theorem for groups in 1963 was 255 pages long, and considered at the time to be very long.  Currently there are a fair number of papers several hundred pages long [the classification of finite simple groups, dated 2004, covers many papers with total length around 10000 to 20000 pages].  The length of papers (or proofs extending over several papers) tends thus to increase with time, in agreement with Turing's result on the halting problem [11].  This explains the growing use of proofs involving computers in an essential way (the 4-colour theorem, the proof of the Kepler conjecture on sphere packings, etc.)
\medskip\noindent
{\bf 5. Understandability.}	
\medskip
	If we look at mathematicians of a certain period (say the early 20-th century: Poincar\'e, Hilbert, G\"odel, von Neumann, etc.) we see that they have very different intellectual profiles.  We can make sense of the statistical difference between scientific ability and other human aptitudes as an effect of evolutionary selection.  In fact scientific ability is a side effect of the development of the human brain in the last one or two million years, followed by the acquisition of language, and quite recently by the ability to count and write.  There is no strong evolutionary pressure in favour of mathematical ability.  This explains why the human ability to produce mathematical arguments understandable by other humans is so variable.  Understandability is a natural human concept which is not easy to formalise but we have seen how to make sense of it in the case the human mathematical structures.
\medskip
	Human short-term memory (working memory) is quite limited (to about 7 objets) compared with the huge short-term memory of computers.  We have argued that this and the use of natural languages explains the introduction of concepts or structures in mathematics.  Mathematical thinking naturally involves combinatorial functioning associated with working memory, and also longer term memory in a way which is not fully clarified.  Understanding a mathematical situation therefore involves a conscious or unconscious functioning of the brain which remains to be described clearly, but we have at least outlined in the previous section how to make sense of the understanding of a mathematical proof.
\medskip
	The human brain is a remarkable information-processing system.  It has its limitations but is very flexible and powerful.  It uses poorly understood methods, including natural languages and concept creation, to deal with complex inter-human communication and reaching very complex results in mathematics.
\medskip
	We have seen how understandable mathematical statements (like the 4-colour theorem) can have a proof that escapes human understanding.  (The statement is understandable -- in the sense of checking that certain rules are satisfied by a mathematical text -- even if the proof is not.)  One can also imagine results in formal language which, because of length and complexity, escape human understanding.  Yet such results might be logically important for the structure of the consequences of the basic mathematical axioms.  The set of theorems properly understandable by humans, and the set of mathematical structures that human mathematicians have introduced appear thus to be a rather limited view of the whole of mathematics.
\medskip\noindent
{\bf 6. Non-mathematical sciences.}
\medskip
	Let us now turn to non-mathematical sciences.  These use observational or experimental protocols to make contact with reality at a certain approximation.  We have thus pieces of approximately observed reality, and a scientific theory consists in identifying a piece of observed reality with a certain logical-mathematical structure.  [This may be a list of biological genera and species, or a system of variables called velocity, acceleration, etc.].
\medskip
	Here are some questions which arise in the special case of physics-astronomy.  Does the standard model of particles account for experiments with strong interactions?  Does general relativity account for the large scale structure of the universe?  Is the dynamics of planets in the solar system chaotic?  It is seen that these questions are of different natures, depending on observational or experimental results and on the numerical study of mathematical models.  Various ingredients entering the study of these questions vary with time in different manners.  One ingredient is political decisions concerning large scale experimental or observational projects, like the decision to end the Superconducting Super Collider project in 1993.  Another ingredient is progress in computer technology; for instance Moore's law that describes the increase of the power of computers with time in the period 1970--2020.  Therefore the rate of progress in non-mathematical sciences is a complicated issue and we lack general statements covering past and present periods.  The same can be said for the kind of mathematical progress which involves computers in an essential way.
\bigskip\noindent
{\bf References.}
\medskip\noindent 
[1]  W. Feit and J.G. Thompson.  ``Solvability of groups of odd order.''  Pacific J. Math. {\bf 13},775-1027(1963).
\medskip\noindent
[2]  J. Hadamard.  {\it The Psychology of Invention in the Mathematical Field.}  Princeton University Press, Princeton, 1945.
\medskip\noindent
[3]  T.C. Hales. ``Formal proof.''  AMS Notices {\bf 55},1370-1380(2008).
\medskip\noindent
[4]  T.C. Hales.  ``Developments in formal proofs'', S\'eminaire Bourbaki n$^o$ 1086(2014).
\medskip\noindent
[5]  G.H. Hardy.  {\it A mathematician's apology.}  Cambridge U. P., Cambridge, 1940.
\medskip\noindent
[6]  A. Jaffe and F. Quinn.  ``Theoretical mathematics: towards a cultural synthesis of mathematics and theoretical physics''. Bull. Amer. Math. Soc. {\bf 29},1-13(1993).
\medskip\noindent
[7]  H. Poincar\'e.  {\it Science et m\'ethode.}  Flammarion, Paris, 1908.
\medskip\noindent
[8]  Y. Pomeau, M. Le Berre and J. Ginibre.  ``Ultimate statistical physics: fluorescence of a single atom.''  arXiv:1605.09253v2 [physics.atom-ph].
\medskip\noindent
[9]  D. Ruelle  ``Conversations on mathematics with a visitor from outer space'' pp. 251-259 in {\it Mathematics: Frontiers and Perspectives.}  Amer. Math. Soc., Providence, RI, 2000.
\medskip\noindent
[10]  M.A.F. Sanj\'uan. ``Artificial intelligence, chaos, prediction and understanding in science.''  I.J. Bifurcation and Chaos {\bf 31},2150173(2021).
\medskip\noindent
[11]  A.M. Turing  ``On computable numbers, with an application to the Entscheidungsproblem'', Proceedings of the London Mathematical Society, s2-42: 230–265; correction ibid., s2-43: 544–546(1937).
\medskip\noindent
[12]  A. Turing.  ``Computing machinery and intelligence.''  Mind {\bf 59},433-460(1950).
\medskip\noindent
[13]  J. von Neumann.  {\it The Computer and the Brain.}  Yale U. P., New Haven, 1958.
\medskip\noindent
[14]  A.N. Whitehead.  {\it Process and reality.} Free Press, New York, 1979 (p. 39)
\end